\documentclass[twocolumn,amsmath,amssymb,aps,secnumarabic,%
    nofootinbib,groupedaddress]{revtex4}
\usepackage{times,mathptmx,amsthm,pstricks,pst-plot}
\newtheorem{theorem}{Theorem}[section]
\newtheorem{lemma}[theorem]{Lemma}
\newtheorem{conjecture}[theorem]{Conjecture}
\newtheorem{corollary}[theorem]{Corollary}
\theoremstyle{definition}
\theoremstyle{remark}
\newtheorem*{remark}{Remark}
\newtheorem*{remarks}{Remarks}

\newcommand{\R}{\mathbb{R}}
\newcommand{\C}{\mathbb{C}}
\newcommand{\eps}{\epsilon}
\renewcommand{\tilde}{\widetilde}

\renewcommand{\tensor}{\otimes}

\newcommand{\ve}{\vec{e}}
\newcommand{\vp}{\vec{p}}
\newcommand{\vv}{\vec{v}}
\newcommand{\vw}{\vec{w}}
\newcommand{\vx}{\vec{x}}
\newcommand{\vy}{\vec{y}}
\newcommand{\tvx}{\tilde{\vec{x}}}
\newcommand{\tvy}{\tilde{\vec{y}}}

\newcommand{\ISO}{\mathrm{ISO}}
\newcommand{\Isom}{\mathrm{Isom}}
\newcommand{\lk}{\mathrm{lk}}
\newcommand{\SO}{\mathrm{SO}}
\newcommand{\so}{\mathrm{so}}
\renewcommand{\bar}{\overline}
\newcommand{\Vol}{\mathrm{Vol}\ }
\renewcommand{\diamond}{\diamondsuit}
\newcommand{\heart}{\heartsuit}

\newcommand{\st}{{\bigm|}}

\newcommand{\ie}{\textit{i.e.}}
\newcommand{\Ie}{\textit{I.e.}}
\newcommand{\eq}[2]{\begin{equation}\label{#1}#2\end{equation}}

\newcommand{\eatline}{\vspace{-\baselineskip}}

\newcommand{\thm}[1]{Theorem~\ref{#1}}
\renewcommand{\sec}[1]{Section~\ref{#1}}
\newcommand{\lem}[1]{Lemma~\ref{#1}}
\newcommand{\cor}[1]{Corollary~\ref{#1}}
\newcommand{\conj}[1]{Conjecture~\ref{#1}}
\newcommand{\fig}[1]{Figure~\ref{#1}}

\newenvironment{fullfigure}[2]
    {\begin{figure}[htb]\begin{center}\def\fullfiga{#1}\def\fullfigb{#2}}
    {\vspace{\baselineskip}\caption{\fullfigb.}\label{\fullfiga}
    \end{center}\end{figure}}

\newgray{gray5}{.5}
\newgray{gray6}{.6}
\newgray{gray7}{.7}
\newgray{gray8}{.8}
\newgray{gray9}{.9}
\psset{hatchwidth=.4pt,hatchsep=3pt}
\psset{linewidth=.4pt,dash=3pt 3pt,doublesep=.05,dotsize=1pt 5,arrowsize=2pt 3}
\SpecialCoor

\begin{document}
\title{From the Mahler conjecture to Gauss linking integrals}

\author{Greg Kuperberg}
\thanks{This material is based upon work supported by the National Science
    Foundation under Grant No. 0606795}
\affiliation{Department of Mathematics, University of
    California, Davis, CA 95616}

\begin{abstract}
\centerline{\textit{\normalsize Dedicated to my father,
    on no particular occasion}}
\vspace{\baselineskip}
We establish a version of the bottleneck conjecture, which in turn implies a
partial solution to the Mahler conjecture on the product $v(K) = (\Vol K)(\Vol
K^\circ)$ of the volume of a symmetric convex body $K \in \R^n$ and its polar
body $K^\circ$.  The Mahler conjecture asserts that the Mahler volume $v(K)$ is
minimized (non-uniquely) when $K$ is an $n$-cube.  The bottleneck conjecture (in
its least general form) asserts that the volume of a certain domain $K^\diamond
\subseteq K \times K^\circ$ is minimized when $K$ is an ellipsoid. It implies
the Mahler conjecture up to a factor of $\bigl(\frac{\pi}4\bigr)^n \gamma_n$,
where $\gamma_n$ is a monotonic factor that begins at $\frac4{\pi}$ and
converges to $\sqrt{2}$.  This strengthens a result of Bourgain and Milman, who
showed that there is a constant $c$ such that the Mahler conjecture is true up
to a factor of $c^n$.

The proof uses a version of the Gauss linking integral to obtain a constant
lower bound on $\Vol K^\diamond$, with equality when $K$ is an ellipsoid. It
applies to a more general conjecture concerning the join of any two necks of
the pseudospheres of an indefinite inner product space.  Because the
calculations are similar, we will also analyze traditional Gauss linking
integrals in the sphere $S^{n-1}$ and in hyperbolic space $H^{n-1}$.
\end{abstract}
\maketitle

\section{Introduction}

If $K \subset \R^n$ is a centrally symmetric convex body, let $K^\circ$ denote
its dual or polar body.  (It is also the unit ball of the norm dual to the one
defined by $K$.)  The product of the volumes
$$v(K) = (\Vol K)(\Vol K^\circ)$$
is known as the Mahler volume of $K$.  Since it is both continuous in $K$ and
affinely invariant, it achieves a finite maximum and a non-zero minimum in each
dimension $n$.  Mahler \cite{Mahler:convexe} conjectured an upper bound and a
lower bound for the value of $v(K)$ in each fixed dimension.  The upper bound
was proven by Santal\'o and is known as the Blaschke-Santal\'o inequality
\cite{Blaschke:vorlesungen2,Santalo:cuerpos,Saint-Raymond:volume}:

\begin{theorem}[Blaschke, Santal\'o] In any fixed dimension $n$, $v(K)$ is
uniquely maximized by ellipsoids. \label{th:santalo}
\end{theorem}

The lower bound is still a conjecture:

\begin{conjecture}[Mahler] In any fixed dimension $n$, $v(K)$ is minimized among
centrally symmetric convex bodies by the cube $C_n = [-1,1]^n$.
\label{c:mahler}
\end{conjecture}

The minimization problem is considered harder because the conjectured minimum
(or the real minimum) is a much more complicated shape than a round sphere.
For instance, Saint-Raymond \cite{Saint-Raymond:volume} observed that the
Mahler volume of an $n$-cube $C_n$ is tied not only by the volume of the dual
$C_n^*$, a cross-polytope, but also by other polytopes when $n$ is large.
(Nonetheless, the Mahler conjecture has been established in some special
cases, including 1-unconditional convex bodies \cite{Saint-Raymond:volume}
and zonoids \cite{Reisner:zonoids,GMR:zonoids}.)

In a noted paper, Bourgain and Milman showed that the Mahler conjecture is
true up to an exponential factor \cite{BM:volume}:

\begin{theorem}[Bourgain, Milman] There is a constant $c>0$ such that
for any $n$ and any centrally-symmetric convex body $K$ of dimension $n$,
    $$v(K) \ge c^n v(C_n).$$ \label{th:bm}
\vspace{-\baselineskip}\end{theorem}

The proof of \thm{th:bm} has been simplified by Pisier \cite{Pisier:volume} and
independently by Milman.  Although all of the proofs technically construct the
constant $c$, no good value for it is currently known.  Note that the Mahler
volume of a sphere and a cube differ only by a factor of $c^n$.  So the
Bourgain-Milman theorem says that all Mahler volumes in $n$ dimensions only
differ by some exponential factor.  (Indeed, the Bourgain-Milman theorem
also holds without central symmetry; see below.)

In this article, we will establish a better value of the constant $c$ by
comparing the Mahler volume with another volume.  We will minimize the volume of
a certain body $K^\diamond \subset K \times K^\circ$.  Define the subsets
$$K^\pm = \{(\vx,\vy) \in K \times K^\circ \st \vx \cdot \vy = \pm 1\}$$
of the hyperboloids
$$H^\pm = \{(\vx,\vy) \st \vx \cdot \vy = \pm 1\}$$
in $\R^n \times \R^n$.  Suppose that $K$ and $K^\circ$ are both positively
curved. (\Ie, they have no boundary points with zero extrinsic curvature in
any direction.)  We previously showed that $K^+$ is a spacelike section of
$H^+$, that $K^-$ is a timelike section of $H^-$, and that their inclusions are
homotopy equivalences \cite{Kuperberg:bottleneck}.  Submanifolds of $H^+$ and
$H^-$ of this type are called \emph{necks}. The body $K^\diamond$ is the filled
join, defined in \sec{s:review}, of the two necks $K^+$ and $K^-$. (In
Reference~\citealp{Kuperberg:bottleneck}, $K^\diamond$ was the convex hull of
$K^+$ and $K^-$ and $K^\heart$ was the filled join.  Since the filled join is
the real object of study, we have renamed it $K^\diamond$.  We do not know if it
is ever non-convex.)

\begin{theorem}[Main theorem] Let $N^+$ and $N^-$ by any two necks of the
positive and negative unit pseudospheres $H^+$ and $H^-$ in any indefinite inner
product space $\R^{(a,b)}$, and let $N^\diamond$ be their filled join.  Then
$\Vol N^\diamond$ is minimized when $N^+$ and $N^-$ are flat, orthogonal, and
centered at the origin.
\label{th:abnecks}
\end{theorem}

Because of the geometry of the necks $N^+$ and $N^-$, we previously called
\thm{th:abnecks} ``the bottleneck conjecture''.  We will interpret
\thm{th:abnecks} as a new isoperimetric-type inequality for convex bodies:

\begin{corollary} If $K \subset \R^n$ is a centrally symmetric convex body, then
$\Vol K^\diamond$ is minimized when $K$ is an ellipsoid.
\label{c:iso}
\end{corollary}

\cor{c:iso} follows immediately from  \thm{th:abnecks} when $K$ and $K^\circ$
are positively curved, because necks $K^+$ and $K^-$ that come from a symmetric
convex body $K$ are flat, orthogonal, and centered when $K$ is an ellipsoid.  By
continuity, the bound must also hold for arbitrary $K$.

On the one hand, $\Vol K^\diamond$ is minimized when the Mahler volume $v(K)$ is
maximized. On the other hand, $K^\diamond$ is only moderately exponentially
smaller than its superset $K \times K^\circ$. Explicitly, let $B_n$ be the round
unit ball in $\R^n$.  Then
$$\Vol B_n^\diamond = \frac{2^n(n!)^2\pi^n}{(2n)!(n/2)!^2}$$
because $B_n$ is the join of two round, orthogonal $n$-balls of radius
$\sqrt{2}$ and
$$\Vol B_n = \frac{\pi^{n/2}}{(n/2)!}.$$
In comparison,
$$v(C_n) = \frac{4^n}{n!}.$$
In conclusion:

\begin{corollary} If $K \subset \R^n$ is a centrally symmetric convex
body, then
$$v(K) \ge \frac{2^n(n!)^2\pi^n}{(2n)!(n/2)!^2}
    = \gamma_n \biggl(\frac{\pi}{4}\biggr)^n v(C_n),$$
where
$$\frac{4}{\pi} \le \gamma_n = \frac{n!^32^n}{(2n)!(n/2)!^2} \to \sqrt{2}.$$
\label{c:bound}
\end{corollary}

Thus, our value of $c$ in \thm{th:bm} is $\pi/4$.  Since the
Bourgain-Milman theorem is also often phrased in terms of $v(B_n)$, the Mahler
volume of a round ball, we can also write \cor{c:bound} as
$$v(K) \ge \frac{2^n(n!)^2}{(2n)!} v(B_n) > 2^{-n}v(B_n).$$
Thus, our value for the Bourgain-Milman constant in this form is $\frac12$.
Mahler's conjecture would imply $\frac{2}{\pi}$.

\cor{c:bound} also improves a previous non-asymptotic estimate which
was inspired by the bottleneck conjecture \cite{Kuperberg:convex}.  There we
established that if $K \subset \R^n$ and $n \ge 4$, then
$$v(K) \ge (\log_2 n)^{-n} v(B_n).$$
\cor{c:bound} is strictly better for all $n \ge 4$.

There is also a version of the Mahler conjecture for asymmetric convex bodies
that contain the origin in their interiors.

\begin{conjecture}[Mahler] In any fixed dimension $n$, $v(K)$ is minimized among
convex bodies $K$, whose interiors contain the origin, by a centered simplex
$\Delta_n$.
\end{conjecture}

Although we will only directly study centrally symmetric convex bodies,
\cor{c:bound} yields a corollary for general convex bodies using
a standard technique.

\begin{corollary} If $K \subset \R^n$ is a convex
body with the origin in its interior, then
$$v(K) \ge \frac{4^n(n!)^4\pi^n}{(2n)!^2(n/2)!^2}
    = \delta_n \biggl(\frac{\pi}{2e}\biggr)^n v(\Delta_n),$$
where
$$\frac{2e}{\pi} \le \delta_n = \frac{8^nn!^6e^n}{(2n)!^2(n/2)!^2(n+1)^{n+1}}
    \to \frac{2\pi}{e}.$$
\end{corollary}

Thus, our value for the Bourgain-Milman constant in the asymmetric case is
$\frac{\pi}{2e}$.  

\begin{remark} Even though the Bourgain-Milman inequality holds in
the asymmetric case, the bound in \cor{c:bound} does not.  In particular,
$$\lim_{n \to \infty} \left(\frac{v(\Delta_n)}{v(C_n)}\right)^{1/n}
    = \frac{e}{4} < \frac{\pi}{4}.$$
The proof fails because when $K$ is not centrally symmetric, only
$K^+$ exists as a subset of $K \times K^\circ$.  The region $K^\diamond$
can be formed but its volume may be larger than $v(K)$.  Nonetheless
\thm{th:abnecks} does say something about asymmetric convex bodies.
We can let $N^+ = K_1^+$ for one convex body $K_1$, and let $N^-$
be a natural isometric image of $K_2^+$ for another convex body $K_2$.
Then the statement is that $\Vol N^\diamond$ is minimized when
$K_1$ and $K_2$ coincide, and are a centered ellipsoid.
\end{remark}

\begin{proof} Given a convex body $K$, then $K-K$ is its \emph{difference body},
by definition the set of differences between points in $K$.  It is a centrally
symmetric convex body.  Rogers and Shephard \cite{RS:difference} showed that
$$\Vol K-K \le \binom{2n}{n} \Vol K,$$
with equality if and only if $K$ is a simplex.  The polar body of $K-K$ is best
described by a relation between the corresponding norms on $\R^n$:
$$||\vx||_{(K-K)^\circ} = ||\vx||_{K^\circ} + ||\vx||_{-K^\circ}.$$
Since
$$\Vol K = \int_{S^{n-1}} \frac{1}{n||\vx||_K^n} d\vx,$$
and since the integrand is convex, Jensen's inequality tells us that
$$\Vol (K-K)^\circ \le 2^{-n}(\Vol K^\circ),$$
with equality if and only if $K$ is centrally symmetric.  Combining these
estimates with \cor{c:bound} establishes the claim. The numerical comparison to
the simplex case uses the elementary volume formula
$$v(\Delta_n) = \frac{(n+1)^{n+1}}{(n!)^2}.$$
\end{proof}

\sec{s:linking} is a digression in which we will compute invariant linking forms
for the sphere $S^{n-1}$ and for hyperbolic space $H^{n-1}$.  The computations
are similar to \lem{l:ode} used in the proof of \thm{th:abnecks}.  We also give
a second derivation in the sphere case using elementary geometry and probability
instead of Stokes' theorem.

\acknowledgments

The author would like to thank Tom Ilmanen, Misha Kapovich, Boaz Klartag, Bruce
Kleiner, George Kozlowski, Peter Michor, Vitali Milman, Andrew Waldron, and his
parents for helpful discussions.  The author would especially like to thank
Semyon Alesker, Dennis DeTurck, and Herman Gluck for helping to check the main
result.  Finally, the author would like to thank David Kazhdan and Joseph
Bernstein for introducing him to the Bourgain-Milman theorem in 1987.

\section{A review}
\label{s:review}

In this section, we review some relevant definitions and results in
Reference~\citealp{Kuperberg:bottleneck} and set up our new arguments.  As in
the introduction, let $K \subset \R^n$ be a centrally symmetric convex body such
that $K$ and $K^\circ$ are both positively curved.  (Note that $K$ is positively
curved if and only if $K^\circ$ is.)

Let $V$ be an (non-degenerate) indefinite inner product space; we will write all
inner products as dot products.  Recall that a non-zero vector $\vx \in V$ is
\emph{spacelike} if $\vx \cdot \vx > 0$, \emph{timelike} if $\vx \cdot \vx < 0$,
and \emph{null} if $\vx \cdot \vx = 0$.  (The definitions of spacelike and
timelike can be switched in some contexts.)  The positive and negative unit
pseudospheres of $V$ are the solutions sets to the equations $\vx \cdot \vx =
\pm 1$.  Every indefinite inner product space is isometric to the standard inner
product space $\R^{(a,b)}$ with inner product
\begin{multline*}
\vx \cdot \vy = x_1y_1 + x_2y_2 + \cdots + x_ay_a \\
    - x_{a+1}y_{a+1} - x_{a+2}y_{a+2} - \cdots - x_{a+b}y_{a+b}
\end{multline*}
with signature $(a,b)$, for some $a$ and $b$.

Recall that a \emph{pseudo-Riemannian} manifold $M$ is a smooth manifold with a
smooth field of non-degenerate inner products, which is called a \emph{metric}.
The metric on $M$ restricts to a metric on every submanifold $N \subset M$;
however, the restriction may or may not be non-degenerate.   A submanifold $N$
of a pseudo-Riemannian manifold $M$ is \emph{spacelike}, \emph{timelike}, or
\emph{null} if all of its non-zero tangent vectors are, respectively, spacelike,
timelike, or null in their respective tangent spaces.

For example, the positive and negative pseudospheres $H^\pm$ of
$\R^{(a,b)}$ are pseudo-Riemannian with signature $(a-1,b)$ and $(a,b-1)$.
They are diffeomorphic to $S^{a-1} \times \R^b$ and $\R^a \times S^{b-1}$.
The diffeomorphisms can be chosen so that every fiber of the first factor is
spacelike and every fiber of the second factor is timelike.

If $A$ and $B$ are two sets in $\R^n$, we define their \emph{geometric join} $A
* B$ to be the union of all line segments that connect a point in $A$ and a
point in $B$.  The geometric join is a kind of partial convex hull.  If $M
\subset \R^n$ is a closed manifold of codimension 1, we define the
\emph{filling} $\bar{M}$ of $M$ to be the compact region in $\R^n$ that it
encloses.  If $A * B$ is such a manifold, then $\bar{A * B}$ is the \emph{filled
join} of $A$ and $B$.

The vector space $\R^n \times \R^n$ has a relevant inner product (or dot
product), with signature $(n,n)$, given by the formula
\eq{e:ip}{(\vx_1,\vy_1)\cdot(\vx_2,\vy_2) =
    \frac{\vx_1\cdot\vy_2+\vx_2\cdot\vy_1}2.}
The hyperboloids $H^+$ and $H^-$ from the introduction are the positive and
negative unit pseudospheres with respect to this inner product.  Note that this
inner product does not have determinant $\pm 1$.  We will use the standard
volume structure on $\R^n \times \R^n$ rather than the one induced by
\ref{e:ip}.  (It may have been a bit clearer to put $K$ in a vector space $V$
and then define this inner product on $V \oplus V^*$.  We will use the notation
$K \subset \R^n$ because it is standard in finite-dimensional convex geometry.)

\begin{fullfigure}{f:diamond}
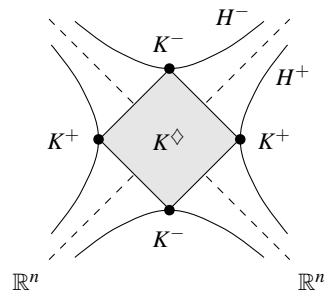
    {The geometry of $K^+$, $K^-$, and $K^\diamond$ in $\R^n \times \R^n$}
\pspicture(-2.1,-2.1)(2.1,1.9)
\psline[linestyle=dashed](-1.6,-1.6)(1.6, 1.6)
\psline[linestyle=dashed](-1.6, 1.6)(1.6,-1.6)
\rput{45}(0,0){
\pscurve( 2, .222)( 1.333, .333)( .666, .666)( .333, 1.333)( .222, 2)
\pscurve(-2, .222)(-1.333, .333)(-.666, .666)(-.333, 1.333)(-.222, 2)
\pscurve( 2,-.222)( 1.333,-.333)( .666,-.666)( .333,-1.333)( .222,-2)
\pscurve(-2,-.222)(-1.333,-.333)(-.666,-.666)(-.333,-1.333)(-.222,-2)
\psframe[fillstyle=solid,fillcolor=gray9](-.666,-.666)(.666,.666)
\qdisk( .666, .666){2pt} \qdisk(-.666, .666){2pt}
\qdisk( .666,-.666){2pt} \qdisk(-.666,-.666){2pt}}
\rput(0,0){$K^\diamond$}
\rput(-1.4,0){$K^+$} \rput(1.4,0){$K^+$}
\rput(0,-1.3){$K^-$} \rput(0,1.3){$K^-$}
\rput(-1.9,-1.9){$\R^n$} \rput(1.9,-1.9){$\R^n$}
\rput(1.65,.85){$H^+$}
\rput(.85,1.65){$H^-$}
\endpspicture
\end{fullfigure}

In this section, we will establish the facts that $K^+$ is a neck of $H^+$ (to
be defined precisely below), that $K^-$ is a neck of $H^-$, and that their
filled join $K^\diamond = \bar{K^+ * K^-}$ is a starlike body that lies between
$H^+$ and $H^-$.  The overall geometry is shown in \fig{f:diamond}.  As the
figure suggests, $K^\diamond$ is sometimes all of $K \times K^\circ$; for
example, when $K = C_n$. (A useful exercise is to work out the geometry of
$K^\diamond$ when $K = C_2$.  In this case, $K^+$ and $K^-$ are non-planar
octagons that visit all 16 vertices of $K \times K^\circ$, which is a 4-cube.)
But in other cases, $K^\diamond$ is significantly smaller than $K \times
K^\circ$.

Every point in $K^+$ is a pair consisting of a point $\vx \in K$ and a dual
vector $\vy \in K^\circ$ that represents a supporting hyperplane of $\vx$.  In
this situation, $\vx \in \partial K$ and $\vy \in \partial K^\circ$.  Moreover,
since $K$ and $K^\circ$ are both smooth, $\vx$ uniquely and smoothly determines
$\vy$ and vice versa. Therefore $K^+$ is diffeomorphic to $\partial K$, and in
particular, is diffeomorphic to the sphere $S^{n-1}$.  Likewise
$K^-$ is as well.

\begin{lemma} $K^+$ is a
spacelike submanifold of $H^+$ and $K^-$ is a timelike submanifold of $H^-$.
\end{lemma}

\begin{proof} (Sketch) The proof is from
Reference \citealp{Kuperberg:bottleneck}.  The convex body $K$ has a unique
osculating ellipsoid $E$ at the point $\vx$.  Because $E$ is an ellipsoid, $E^+$
is a flat ellipsoid surface in $H^+ \subset \R^n \times \R^n$, and in particular
is a spacelike manifold.  Since $E$ osculates $K$ at $\vx$, it follows that
$$T_{(\vx,\vy)}K^+ = T_{(\vx,\vy)}E^+.$$
Thus, $K^+$ is spacelike, and likewise $K^-$ is timelike.
\end{proof}

In addition, $K^+$ is a spacelike section of $H^+$, in two senses. Since $H^+$
is a pseudo-Riemannian manifold of signature $(n-1,n)$, $K^+$ is a spacelike
section in the sense of having maximal dimension. But also, the normal
exponential map from any flat section $K^+$ of $H^+$ (such as $E^+$ if $E$ is an
ellipsoid) makes $H^+$ a bundle of hyperbolic spaces over a sphere.  Then $K^+$
is a section of this bundle, because it is both spacelike and isotopic to
$K^+$.  ($K^+$ is isotopic to $E^+$ by deforming $K$ to $E$, and then
all flat sections are isotopic.)  Likewise $K^-$ has the same properties
in $H^-$.

If $V$ is any indefinite inner product space and $H^\pm$ are its positive and
negative unit pseudospheres, define a \emph{neck} of $H^\pm$ to be a spacelike
or timelike submanifold $N^\pm$ which is a section in both of the above senses.
If $V$ has signature $(a,b)$, then $N^+$ has dimension $a-1$ and $N^-$ has
dimension $b-1$.

For the main result, we will let $V = \R^{(a,b)}$.  Each spacelike $a$-plane in
$\R^{(a,b)}$ has a canonical orientation, as does each timelike $b$-plane.  The
group of isometries of $\R^{(a,b)}$ that preserve its total orientation is
$\SO(a,b)$.  The subgroup that separately preserves the orientations of maximal
spacelike and timelike plane is $\ISO(a,b)$; it is a connected Lie group.

Let $N^+$ and $N^-$ be necks of $H^+$ and $H^-$.

\begin{lemma} The geometric join $N^+ * N^-$ is the boundary of a starlike
body in $\R^n \times \R^n$.
\label{l:starlike} \end{lemma}

\begin{proof} The argument is a refinement of the one in
Reference~\citealp{Kuperberg:bottleneck}, where it was argued that $K^+ * K^-$
is starlike except possibly on a set of measure 0.

Let $N^+ \circledast N^-$ denote the abstract join of $N^+$ and $N^-$ as defined
in topology (and usually written $N^+ * N^-$), and let
$$\Phi: N^+ \circledast N^- \to \R^{(a,b)}$$
be the obvious continuous map with image $N^+ * N^-$, given by the formula
$$\Phi(\vx,\vy,t) = (1-t)\vx + t\vy.$$
Let
$$\rho:\R^{(a,b)} \to S^{a+b-1}$$
be the tautological radial projection.

The composition $\rho \circ \Phi$ is a map between $(a+b-1)$-spheres.  We claim
that it has degree 1, that it is smooth with non-zero
Jacobian except on $N^+$ and $N^-$, and that it is a local homeomorphism at
$N^+$ and $N^-$.  The last two claims imply that $\rho \circ \Phi$ is a covering
map and therefore a bijection by the first claim. This is equivalent to the
assertion of the lemma.

It is elementary that $\rho \circ \Phi$ is smooth on the smooth points of its
domain.

The claim that $\rho \circ \Phi$ has degree 1 is not really necessary, because
$S^{a+b-1}$ is simply connected unless $a=b=1$, in which case it is easy to show
that $N^+ * N^-$ is a starlike quadrilateral.   One way to show it is that $\rho
\circ \Phi$ varies continuously as $N^+$ and $N^-$ are varied.  If we take $N^+$
and $N^-$ to be orthogonal, flat, and centered, then it can be seen directly
that $N^+ * N^-$ is starlike and indeed convex, which implies that $\rho \circ
\Phi$ has degree 1.

To confirm that $\rho \circ \Phi$ has non-zero Jacobian away from $N^+$
and $N^-$, let $\vx \in N^+$ and $\vy \in N^-$, and let
$$\vv_1,\vv_2,\ldots,\vv_{a-1} \in T_{\vx}N^+
\qquad \vw_1,\vw_2,\ldots,\vw_{b-1} \in T_{\vy}N^-$$
be bases. Suppose that they are chosen so that
$$\vx,\vv_1,\vv_2,\ldots,\vv_{a-1}
    \qquad \vy,\vw_1,\vw_2,\ldots,\vw_{b-1}$$
are both positively oriented bases of maximal timelike and spacelike planes.
Given coordinates on $N^+$ and $N^-$ that induce these tangent bases, the
Jacobian of $\rho \circ \Phi$ at the point $(\vx,\vy,t)$ is
\begin{multline*}
J(\vx,\vy,t) = \\
\frac{\det \bigl[
    (1-t)\vv_1 , \ldots , (1-t)\vv_{a-1} ,
    t\vw_1 , \ldots , t\vw_{b-1} , \vy-\vx , (1-t)\vy + t\vx \bigr]}
    {||(1-t)\vv + t\vw||_2^{a+b}} \\
    = \frac{(1-t)^{a-1}t^{b-1}
    \det \bigl[\vv_1,\ldots,\vv_{a-1},\vw_1,\ldots,\vw_{b-1},\vy,\vx\bigr]}
    {||(1-t)\vv + t\vw||_2^{a+b}},
\end{multline*}
and so does not vanish.

\begin{fullfigure}{f:hole}
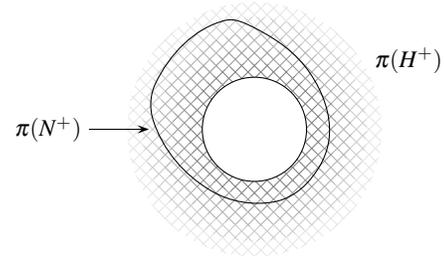
    {$\pi(N^+)$ is starlike and encloses the hole of $\pi(H^+)$}
\pspicture(-3,-1.7)(3,1.7)
\pscircle[linestyle=none,fillstyle=crosshatch*,hatchcolor=gray9](0,0){1.7}
\pscircle[linestyle=none,fillstyle=crosshatch*,hatchcolor=gray8](0,0){1.5}
\pscircle[linestyle=none,fillstyle=crosshatch*,hatchcolor=gray7](0,0){1.3}
\pscircle[linestyle=none,fillstyle=crosshatch*,hatchcolor=gray6](0,0){1.1}
\pscircle[linestyle=none,fillstyle=crosshatch*,hatchcolor=gray5](0,0){.9}
\pscircle[fillstyle=solid,fillcolor=white](0,0){.7}

\psccurve(1.05;45)(1.35;90)(1.5;105)(1.45;155)(1.35;175)(1;235)(1;285)
    (1;315)(1;355)
\rput[l](1.6,.92){$\pi(H^+)$}
\rput[r](-2.3,.0){$\pi(N^+)$}
\psline{->}(-2.2,0)(-1.4,0)
\endpspicture
\end{fullfigure}

Finally at a point $\vy \in N^-$, the map $\Phi$ has a local conelike
structure.  It is (tautologically) smooth and non-singular along the tangent
space $T_{\vy}N^-$, while the direction of $\vy$ is annihilated by the radial
projection $\rho$.  In the remaining directions, the map $\Phi$ is  locally a
conical extension of the restriction of the projection
$$\pi:\R^{(a,b)} \to \bigl(T_{\vy}N^- + \langle\vy\rangle\bigr)^{\perp}$$
to the other neck $N^+$.  Because $N^+$ is a spacelike section, $\pi(N^+)$ must
be starlike.  This implies that its cone extension is a homeomorphism, and
therefore that $\rho \circ \Phi$ is a local homeomorphism at $\vy$.  The same
argument applies to $\vx \in N^+$.
\end{proof}

\lem{l:starlike} implies that $N^+$ and $N^-$ admit a filled join
$$N^\diamond = \bar{N^+ * N^-}.$$
The lemma also means that the volume of $N^\diamond$ can be
expressed as an integral, by the method of infinitesimals:
\eq{e:int}{\Vol N^\diamond = \frac{(a-1)!(b-1)!}{(a+b)!}
    \int_{(\vx,\vy) \in N^+ \times N^-} \hspace{-3.5em} \vx \wedge \vy \wedge
    d\vx^{\wedge a-1} \wedge d\vy^{\wedge b-1}.}

\begin{remark} Equation~\eqref{e:int} is different in several ways from the
formula in Reference~\citealp{Kuperberg:bottleneck}.  The most important change
is that the constant factor was not correct previously. In addition, the
previous formula was for the volume of $K^\heart$, which is here called
$K^\diamond$.  Finally the old formula had $d\vx$ and $d\vy$ instead of their
wedge powers.  While this was plausible notation, the wedge power notation is
more literally correct and will be important to the proof of \thm{th:abnecks}.
\end{remark}

Equation~\eqref{e:int} is technically only true up to sign, but we can make it
exactly true with suitable orientations for the necks.

The differential forms in equation~\eqref{e:int} are vector-valued and the wedge
products are really ``double wedges'' in the algebra
$$A = \Lambda^*(\R^{(a,b)}) \tensor \Omega^*(\R^{(a,b)}).$$
Even though both factors are graded-commutative algebras, we define
multiplication in $A$ by tensoring them as ungraded algebras. If
$\omega_1,\omega_2 \in A$ have degrees $(j_1,k_1)$ and $(j_2,k_2)$, then
$$\omega_1 \wedge \omega_2 = (-1)^{j_1j_2 + k_1k_2} \omega_2 \wedge \omega_1.$$
For example,
$$\vx \wedge \vy = -\vy \wedge \vx \qquad
    \vx \wedge d\vy = -d\vy \wedge \vx \qquad
    d\vx \wedge d\vy = d\vy \wedge d\vx.$$
Equation~\eqref{e:int} is established by assuming simplex shapes for the
infinitesimal elements $d\vx$ and $d\vy$ of $N^+$ and $N^-$, so that the
corresponding infinitesimal element of $N^\diamond$ is also a simplex (subtended
by $\vx$, $d\vx$, $\vy$, and $d\vy$).  Strictly speaking, the integral on the
right side of \eqref{e:int} is not scalar, but rather has a value in
$\Lambda^{a+b}(\R^{(a,b)})$.  However, it can be read as a scalar because
$\Lambda^{a+b}(\R^{(a,b)})$ is 1-dimensional; we endow it with the standard
basis vector
$$\ve_1 \wedge \ve_2 \wedge \cdots \wedge \ve_{a+b}.$$

\section{Proof of \thm{th:abnecks}}
\label{s:proof}

The idea of the proof is to find a topological lower bound for the integral in
equation \eqref{e:int}.  We can compactify the union $H^+ \cup H^-$ by radially
projecting it onto the Euclidean-unit sphere $S^{a+b-1}$ in $\R^{(a,b)}$,
similar to the proof of \lem{l:starlike}.  The image is dense; the points in
$S^{a+b-1}$ that are missed are null vectors.  Thus we have a compactification
$$\bar{H^+ \cup H^-} \cong S^{a+b-1}.$$
Our lower bound will be proportional to a Gauss-type integral for the linking
number of $N^+$ and $N^-$ in $\bar{H^+ \cup H^-}$, which is necessarily 1.  Our
Gauss linking form will be invariant under the non-compact group $\SO(a,b)$
rather than the usual rotation group $\SO(a+b)$.

If $X$ and $Y$ are two manifolds, then differentials on $X \times Y$ are doubly
graded according to their degrees in the $X$ and $Y$ directions.  The full
differential $d$ also splits as $d = d_x + d_y$ according to differentiation
along $X$ and $Y$ separately.  This refines the de Rham complex
$\Omega^*(X \times Y)$ into a double complex.  In our case, $X = H^+$ and $Y =
H^-$.  As before, we will use the vector coordinates $\vx \in H^+$ and $\vy \in
H^-$ using the inclusions of $H^\pm$ in $\R^{(a,b)}$.

\begin{lemma} Let $\omega$ be an $\ISO(a,b)$-invariant differential form on $H^+
\times H^-$ of degree $(a-1,b-1)$, and let $N^+$ and $N^-$ be closed, smooth
submanifolds or cycles of dimension $a$ and $b$ in $H^+$ and $H^-$. Also suppose
that $a,b>1$.  Then
$$\int_{N^+ \times N^-} \omega$$
is invariant under chain homotopy of $N^+$ and $N^-$ if and only if
$d_yd_x \omega = 0$.
\label{l:dydx} \end{lemma}

\begin{remarks}
Since $d_x$ and $d_y$ anticommute, the condition on $\omega$ is symmetric in
$\vx$ and $\vy$.  In the context of double complexes, $\omega$ can be called
``weakly closed''.

Our form $\omega$ will actually be $\SO(a,b)$-invariant.  But $\SO(a,b)$ is not
connected when $a,b>0$.  Its connected subgroup $\ISO(a,b)$ is more natural for
this lemma (and not really different).

DeTurck and Gluck \cite{DG:personal} show that the double integral of a similar
form $\omega$ (see \sec{s:linking}) is chain-homotopy invariant if and only if
it is weakly closed and
$$d_x\omega = d_y \sigma \qquad d_y\omega = d_x \kappa$$
for some forms $\sigma$ and $\kappa$.  This criterion does not
require $\omega$ to be $\ISO(a,b)$-invariant.
\end{remarks}

\begin{proof} We will first consider the ``if'' direction, which is the
one that we will need.  Let
$$\omega_x = \int_{N^-} \omega.$$
Then
$$d\omega_x = \int_{N^-} d_x\omega$$
is a differential form on $H^+$ of degree $a$.  By Stokes' theorem, it is
invariant under chain homotopy of $N^-$, because $d_y$ annihilates the
integrand.  In particular, $d\omega_x$ is $\ISO(a,b)$-invariant.

But there are no non-zero $\ISO(a,b)$-invariant $a$-forms on $H^+$.  They are
sections of $\Lambda^a(T^*H^+)$.  Each fiber of $T^*H^+$ is a representation of
the point stabilizer $\ISO(a-1,b)$ and is isomorphic to the defining
representation $V$. The exterior power $\Lambda^a(V)$ has no invariant vectors
unless $\dim V = a$. This can be checked by passing to the action of the
complexified Lie algebra $\so(a+b-1,\C)$ on $\Lambda^a(V) \tensor \C$.  The
claim is then a basic calculation in the representation theory of complex simple
Lie algebras. Namely, $\Lambda^a(V) \tensor \C$ is an irreducible representation
of $\so(a+b-1,\C)$ unless $a = b-1$, in which case it is a direct sum of two
irreducible representations.  The only case with $a,b > 0$ in which it is the
trivial representation is when $b = 1$, which we have excluded by hypothesis.

Thus, $\omega_x$ is a closed form on $H^+$, so
$$\int_{N^+} \omega_x = \int_{N^+ \times N^-} \omega$$
is invariant under chain homotopy of $N^+$.  By the same argument, it is also
invariant under chain homotopy of $N^-$.  This concludes the ``if'' part
of the lemma.

To sketch the ``only if'' part, suppose that $N_0^+$ and $N_1^+$ are two necks
of $H^+$ that agree except for a locus which is the boundary of an
$a$-dimensional blister $L^+$.  (By a \emph{blister}, we mean any manifold with
cusp boundary that connects two smooth manifolds that are only partly
disjoint.)  Likewise, suppose that $N_0^-$ and $N_1^-$ are two necks of $H^-$
that differ only at a $b$-dimensional blister $L^-$.  By two applications of
Stokes' theorem,
$$0 = \sum_{p,q \in {0,1}} (-1)^{p+q} \int_{N_p^+ \times N_q^-} \omega
    = \int_{L^+ \times L^-} d_yd_x \omega.$$
But then, the blisters $L^+$ and $L^-$ can be arbitrarily small and
can lie in any direction, so that in the limit,
$$d_yd_x \omega = 0.$$
\end{proof}

\begin{lemma} If
$$\omega = \phi(\vx \cdot \vy)\vx \wedge \vy \wedge d\vx^{\wedge a-1} \wedge
    d\vy^{\wedge b-1} $$
is an $\SO(a,b)$-invariant differential form on $H^+ \times H^-$ in
$\R^{(a,b)}$, then $d_yd_x\omega = 0$ if and only if $f(\alpha) =
\phi(\sinh(\alpha))$ satisfies the ordinary differential equation
$$f'' + (a+b)(\tanh \alpha)f' + abf = 0.$$
\label{l:ode} \eatline \end{lemma}

\begin{proof}
We compute:
\begin{multline*}
d_x\omega =  (\vy \cdot d\vx)\phi'(\vx \cdot \vy)
    \vx \wedge \vy \wedge d\vx^{\wedge a-1} \wedge d\vy^{\wedge b-1} \\
    - \phi(\vx \cdot \vy) \vy \wedge d\vx^{\wedge a} \wedge d\vy^{\wedge b-1}.
\end{multline*}
Then
\begin{align}
d_yd_x\omega = & (\vx \cdot d\vy)(\vy \cdot d\vx) \phi'' \vx \wedge \vy
    \wedge d\vx^{\wedge a-1} \wedge d\vy^{\wedge b-1} \nonumber \\
    & + (d\vy \cdot d\vx)\phi' \vx \wedge \vy \wedge d\vx^{\wedge a-1} \wedge
    d\vy^{\wedge b-1} \nonumber \\
    & - (\vy \cdot d\vx)\phi' \vx \wedge d\vx^{\wedge a-1}
    \wedge d\vy^{\wedge b} \nonumber \\
    & - (\vx \cdot d\vy)\phi' \vy \wedge d\vx^{\wedge a}
    \wedge d\vy^{\wedge b-1} \nonumber \\
    & - \phi d\vx^{\wedge a} \wedge d\vy^{\wedge b} = 0.
\label{e:ddomega}
\end{align}
In addition, the constraints
$$\vx \cdot \vx = 1 \qquad \vy \cdot \vy = -1$$
yield the differential conditions
\eq{e:diffcond}{\vx \cdot d\vx = 0 \qquad \vy \cdot d\vy = 0.}

The form $\omega$ is explicitly invariant under the diagonal action of
$\SO(a,b)$, which is transitive on pairs of vectors $(\vx,\vy) \in H^+ \times
H^-$ with a prescribed dot product.  (Indeed, the given form of $\omega$ is the
general expression for an invariant form of degree $(a-1,b-1)$.)  Therefore we
can check the condition $d_yd_x\omega = 0$ for the standard vectors
$$\vx = (1,0,0,\ldots,0) \qquad
    \vy = (\sinh \alpha,0,0,\ldots,\cosh \alpha,\ldots,0).$$
Here the two non-zero coordinates of $\vy$ are $y_1$ and $y_{a+1}$ and $\alpha$
is a hyperbolic angle (or rapidity). Then the conditions~\eqref{e:diffcond}
simplify to
$$dx_1 = 0 \qquad (\sinh \alpha)dy_1 - (\cosh \alpha)dy_{a+1} = 0.$$
Let $d\tvx$ and $d\tvy$ be $d\vx$ and $d\vy$ with the first and $(a+1)$st
coordinates deleted.  Then equation~\eqref{e:ddomega} becomes (term by term)
\begin{multline*}
(-1)^a d_yd_x\omega = \\
\begin{aligned}
    & (\cosh \alpha)^2\phi'' dy_1 \wedge dx_{a+1} \wedge
    d\tvx^{\wedge a-1} \wedge d\tvy^{\wedge b-1} \\
    & + (\cosh \alpha)\phi' dy_{a+1} \wedge dx_{a+1} \wedge
    d\tvx^{\wedge a-1} \wedge d\tvy^{\wedge b-1} \\
    & - b(\cosh \alpha)\phi' dx_{a+1} \wedge  dy_{a+1} \wedge
    d\tvx^{\wedge a-1} \wedge d\tvy^{\wedge b-1} \\
    & + a(\sinh \alpha)\phi' dy_1 \wedge  dx_{a+1} \wedge
    d\tvx^{\wedge a-1} \wedge d\tvy^{\wedge b-1} \\
    & + ab \phi d\tvx^{\wedge a-1} \wedge d\tvy^{\wedge b-1} = 0
\end{aligned}
\end{multline*}
This expression so far uses the relations
$$dx_1 = 0 \qquad dy_1 \wedge dy_{a+1} = 0$$
to eliminate terms in the expansion.  If we apply the remaining relation
$$(\sinh \alpha)dy_1 - (\cosh \alpha)dy_{a+1} = 0$$
and collect differential factors, we conclude that all
of the terms are proportional to
$$dy_1 \wedge dx_{a+1} \wedge d\tvx^{\wedge a-1} \wedge d\tvy^{\wedge b-1},$$
and the scalar factor yields the equation
$$(\cosh \alpha)^2 \phi'' + (a+b+1)(\sinh \alpha)\phi' + ab\phi = 0.$$
However, the notation of this differential equation is misleading, because
$\phi$ is not directly a function of $\alpha$.  Rather, $\phi = \phi(t)$,
where
$$t = \vx \cdot \vy = \sinh \alpha.$$
If we let $f(\alpha) = \phi(\sinh \alpha)$, then
$$f'' + (a+b)(\tanh \alpha)f' + abf = 0$$
is the corresponding ODE for $f$.
\end{proof}

\begin{fullfigure}{f:ode}
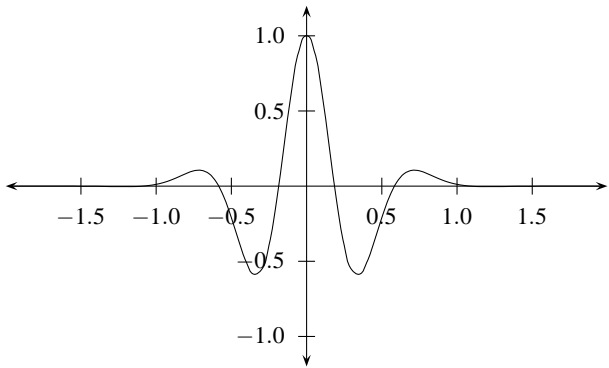
    {The even solution $f(\alpha)$ with $a=b=9$}
\pspicture(-4,-2.5)(4,2.5)
\psset{unit=2cm}
\psaxes[arrows=<->,Ox=0.0,Dx=0.5,Oy=0.0,Dy=0.5](0,0)(-2,-1.2)(2,1.2)
\pscurve
    (-1.95, 0.000)(-1.90,-0.000)(-1.85,-0.000)(-1.80,-0.000)(-1.75,-0.000)
    (-1.70,-0.000)(-1.65,-0.000)(-1.60,-0.000)(-1.55,-0.000)(-1.50,-0.001)
    (-1.45,-0.001)(-1.40,-0.001)(-1.35,-0.002)(-1.30,-0.003)(-1.25,-0.003)
    (-1.20,-0.003)(-1.15,-0.003)(-1.10,-0.001)(-1.05, 0.004)(-1.00, 0.012)
    (-0.95, 0.024)(-0.90, 0.041)(-0.85, 0.062)(-0.80, 0.084)(-0.75, 0.102)
    (-0.70, 0.105)(-0.65, 0.084)(-0.60, 0.027)(-0.55,-0.072)(-0.50,-0.210)
    (-0.45,-0.367)(-0.40,-0.507)(-0.35,-0.586)(-0.30,-0.556)(-0.25,-0.393)
    (-0.20,-0.102)(-0.15, 0.267)(-0.10, 0.633)(-0.05, 0.901)( 0.00, 1.000)
    ( 0.05, 0.901)( 0.10, 0.633)( 0.15, 0.267)( 0.20,-0.102)( 0.25,-0.393)
    ( 0.30,-0.556)( 0.35,-0.586)( 0.40,-0.507)( 0.45,-0.367)( 0.50,-0.210)
    ( 0.55,-0.072)( 0.60, 0.027)( 0.65, 0.084)( 0.70, 0.105)( 0.75, 0.102)
    ( 0.80, 0.084)( 0.85, 0.062)( 0.90, 0.041)( 0.95, 0.024)( 1.00, 0.012)
    ( 1.05, 0.004)( 1.10,-0.001)( 1.15,-0.003)( 1.20,-0.003)( 1.25,-0.003)
    ( 1.30,-0.003)( 1.35,-0.002)( 1.40,-0.001)( 1.45,-0.001)( 1.50,-0.001)
    ( 1.55,-0.000)( 1.60,-0.000)( 1.65,-0.000)( 1.70,-0.000)( 1.75,-0.000)
    ( 1.80,-0.000)( 1.85,-0.000)( 1.90,-0.000)( 1.95, 0.000)
\endpspicture
\end{fullfigure}

\begin{lemma} If $f(\alpha)$ satisfies the differential equation
$$f'' + (a+b)(\tanh \alpha)f' + abf = 0$$
and $f'(0) = 0$, then $|f(\alpha)| < |f(0)|$ when $\alpha \ne 0$.
\label{l:energy} \end{lemma}

\fig{f:ode} shows an example.

\begin{proof} The differential equation for $f$ is that of a damped harmonic
oscillator, where the damping is forwards in time for positive time $\alpha >
0$, and backwards in time for negative time $\alpha < 0$.  Therefore if $f'(0)$,
the oscillator will lose energy in both directions and never again reach $\pm
f(0)$.

In detail, let
$$E = \frac{(f')^2 + abf^2}{2}$$
be the energy of the oscillator.  Then
$$E' = f''f' + abf'f = -(a+b)(\tanh \alpha)(f')^2.$$
Thus $E' \le 0$ when $\alpha > 0$ and vice-versa, with equality only when $f' =
0$.  Thus, $E(\alpha) < E(0)$ for $\alpha \ne 0$.  Moreover,
$$|f| \le \sqrt{\frac{2E}{ab}},$$
with equality when $f' = 0$, as is the case when $\alpha = 0$. These
relations together establish the lemma.
\end{proof}

\begin{proof}[Proof of \thm{th:abnecks}] When $a,b>1$, we combine the
lemmas to establish the theorem.  We want to find the necks $N^+$ and $N^-$
that minimize the integral
$$w(N^+,N^-) = \int_{N^+ \times N^-} \vx \wedge \vy \wedge d\vx^{\wedge a-1}
    \wedge d\vy^{\wedge b-1}.$$
We orient the necks so that the integrand --- not just the integral --- is
positive.  This is possible by the geometry established in the proof of
\lem{l:starlike};  namely, the vector spaces $T_{\vx}N^+ + \langle\vx\rangle$
and $T_{\vy}N^- + \langle\vy\rangle$  are maximal timelike and spacelike
subspaces of $\R^{(a,b)}$.

Let $f$ be the solution to the differential equation
$$f'' + (a+b)(\tanh \alpha)f' + abf = 0$$
with $f(0) = 1$ and $f'(0) = 0$.  Lemmas~\ref{l:dydx}
and \ref{l:ode} say that
$$\ell(N^+,N^-) = \int_{N^+ \times N^-} f(\sinh^{-1} \vx \cdot \vy)
    \vx \wedge \vy \wedge d\vx^{\wedge a-1} \wedge d\vy^{\wedge b-1}$$
is invariant under homotopy of the necks, and is therefore the same for all
necks.  By \lem{l:energy}, $f(\alpha) < f(0) = 1$ when $\alpha \ne 0$.
Therefore
$$\ell(N^+,N^-) \le w(N^+,N^-).$$
Equality is achieved if and only if $\vx \cdot \vy = 0$ always; in other words,
when $N^+$ and $N^-$ are orthogonal, flat, centered necks, as desired.

\lem{l:dydx} is false in the stated generality when $a=1$ or $b=1$. Suppose that
$b=1$ (say). In this case, $H^-$ is a two-sheeted hyperboloid and a valid neck
$N^-$ consists of one point on each sheet.  \lem{l:dydx} does hold, by a
modified argument, if we either require $N^-$ to be a neck, or if we assume an
even solution $f$ to the ODE in \lem{l:ode}.  However, we will give the old
geometric proof from Reference \citealp{Kuperberg:bottleneck} and explain how it
is actually the same proof.  Let $\vy$ be the difference between the two points
of $N^-$; it is necessarily a timelike vector.  After a Lorentz transformation,
$$\vy = (0,0,0,\ldots,0,y)$$ with $y \ge 2$, with equality when $N^-$ is
centered.  Then $v(N^+,N^-)$ is proportional to the area enclosed by the
projection $\pi$ of $N^+$ onto the first $a$ coordinates.  Meanwhile $H^+$
projects onto the outside of the unit sphere $S^{a-1}$ in $\R^a$, and $N^+$ must
also wind around the hole, as in \fig{f:hole}.  The enclosed area is minimized
when $N^+$ is the boundary of the hole, \ie, when $N^+$ is a flat, centered neck
orthogonal to $\vy$. The integral $v(N^+,N^-)$ also has a factor of $y$, and is
minimized when $y=2$, which occurs when $N^-$ is centered.

This completes the proof in all cases.  We remark that if $\rho$ is the radial
projection of $\R^a$ onto $S^{a-1}$ and $\mu$ is Haar measure on $S^{a-1}$, then
the pullback $\pi^*(\rho^*(\mu))$ is the same form $\omega$ that is defined in
Lemmas~\ref{l:dydx} and \ref{l:ode}.  The scalar kernel is
$$f(\alpha) = \frac{1}{(\cosh \alpha)^a},$$
which explicitly satisfies the equation and estimate of \lem{l:energy}.
\end{proof}

\section{Other linking forms}
\label{s:linking}

The calculation of \lem{l:ode} can be adapted to compute Gauss linking forms in
the round $(n-1)$-sphere $S^{n-1}$, or in  hyperbolic space $H^{n-1}$.  If $M$
is one of these two geometries, the form $\omega$ is now an element of
$\Omega^{(a-1,b-1)}(M^{\times 2} \setminus \Delta)$, where $\Delta$ is the
diagonal of $M^{\times 2}$ and $a+b = n$.  The goal is to choose $\omega$ so
that  it is invariant under $\Isom(M)$ and so that its double integral is the
linking number between $(a-1)$-dimensional and $(b-1)$-dimensional submanifolds
$N_1$ and $N_2$ in $M$,
$$\lk(N_1,N_2) = \int_{N_1 \times N_2} \omega,$$
or a proportionality.

This goal is closely related to the problem of finding the propagator of
generalized electromagnetism, or the generalized Biot-Savart law, in curved
spaces.  This question was recently considered by DeTurck and Gluck in $S^3$ and
$H^3$ \cite{DG:electro}. It is not exactly the same question as ours, because if
$\omega$ is such a propagator, then it satisfies a Laplace-type equation instead
of an exterior derivative equation and a symmetry condition.  However, the
propagator has to be an invariant linking form by Ampere's law.

The configuration space $M^{\times 2} \setminus \Delta$ is no longer a Cartesian
product, but it is foliated both horizontally and vertically, so the de Rham
complex $\Omega^{*}(M^{\times 2} \setminus \Delta)$ still refines to a double
complex.  \lem{l:dydx} still holds, except that the chain homotopies on
$N_1$ and $N_2$ must be disjoint.  This allows the integral of  $\omega$ to be
proportional to $\lk(N_1,N_2)$, provided that $d_y d_x \omega = 0$.  No other
invariants are available for the allowed chain homotopies.

We can suppose that $M = S^{n-1}$ is the unit sphere in $\R^{(n,0)}$, or
that $M = H^{n-1}$ is the positive unit pseudosphere in $\R^{(n-1,1)}$.  In
either case, we can write
$$\omega = \phi(\vx \cdot \vy)\vx \wedge \vy \wedge d\vx^{\wedge a-1} \wedge
    d\vy^{\wedge b-1} $$
and perform a calculation similar to the one in \lem{l:ode}.
The calculations are similar enough that we will just give the conclusions.

If $M = S^{n-1}$,
then we let
\begin{align*}
\vx &= (1,0,0,\ldots,0) & \vy &= (\cos \alpha,\sin \alpha,0,0,\ldots,0) \\
f(\alpha) &= \phi(\cos \alpha).
\end{align*}
Then the resulting ODE in $f$ is:
\eq{e:sn1}{f'' + (a+b)(\cot \alpha)f' - abf = 0.}
In this case, the equation has the boundary condition that  $f$ is non-singular
at $\pi$.  With this boundary condition, we obtain
$$\int_{N_1 \times N_2} \omega = f\bigl(\frac{\pi}2\bigr)
(\Vol S^{a-1})(\Vol S^{b-1})\lk(N_1,N_2),$$
where $\Vol S^n$ denotes the volume of the round unit $n$-sphere.
The constant factor is obtained by taking $N_1$ and $N_2$ to be
orthogonal great spheres in $S^{n-1}$ as a test case.  Explicitly,
if $a=b=2$,
\eq{e:s3}{f(\alpha) = \frac{(\pi-\alpha)(\cos \alpha) + (\sin \alpha)}
    {4\pi^2(\sin \alpha^3)}.}
This solution is equivalent to the one given by DeTurck and Gluck
\cite{DG:electro}.

If $M = H^{n-1}$, then we let
\begin{align*}
\vx &= (1,0,0,\ldots,0) & \vy &= (\cosh \alpha,\sinh \alpha,0,0,\ldots,0) \\
f(\alpha) &= \phi(\cosh \alpha).
\end{align*}
Then the resulting ODE in $f$ is:
$$f'' + (a+b)(\coth \alpha)f' + abf = 0.$$
There are two solutions, one even $\alpha$ and the other odd in $\alpha$.
When $a=b=2$, the even and odd solutions are:
$$f(\alpha) = \frac{\cosh \alpha}{(\sinh \alpha)^3} \qquad
    f(\alpha) = \frac{(\sinh \alpha)-\alpha(\cosh \alpha)}{(\sinh \alpha)^3}.$$
If we rescale the odd solution to match the Gauss formula in flat $\R^3$ as
$\alpha \to 0$, we
obtain
$$f(\alpha) = \frac{\cosh \alpha}{4\pi(\sinh \alpha)^3}.$$
This solution is also equivalent to the one given by Gluck and DeTurck.
(It is the analytic continuation of the linking form
on the spherical manifold $\R P^3$.  We could also analytically continue
equation \eqref{e:s3} to obtain a different solution.)

There is also a geometric argument for equation~\eqref{e:s3}.  Let $N_1$ and
$N_2$ be two knots in $S^3$, and let $\vp \in S^3$ be a point that does not lie
on $N_1$, $N_2$, or $-N_1$.  The cone $C = C(N_1,\vp)$ over $N_1$ with vertex
$\vp$ is a well-defined geometric chain in $S^3$, using geodesic segments to
connect $\vp$ to points in $N_1$. Thus the linking number $\lk(N_1,N_2)$ equals
the homological intersection $C(N_1,\vp) \cap N_2$ of $C$ with $N_2$.  This is
still so if we choose $\vp$ at random on $S^3$ with respect to Haar measure.
Therefore if we divide $N_1$ and $N_2$ into infinitesimal segments $d\vx$ and
$d\vy$, the expected value of the homological intersection,
\eq{e:expect}{\omega = E_{\vp}\bigl[C(d\vx,\vp) \cap d\vy\bigr],}
is an invariant linking form.

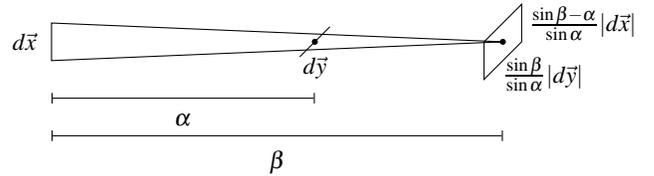
\begin{fullfigure}{f:triangle}{The window of a geometric intersection in $S^3$}
\pspicture(-3.6,-1.5)(4.8,1)
\pspolygon(-3,-.25)(-3,.25)(3,0)
\rput[r](-3.2,0){$d\vx$}
\psline(.7,.2)(.3,-.2)
\rput[t](.5,-.2){$d\vy$}
\psline[arrows=|-|](-3,-.75)(.5,-.75)
\rput[t](-1.25,-.95){$\alpha$}
\psline[arrows=|-|](-3,-1.25)(3,-1.25)
\rput[t](0,-1.45){$\beta$}
\pspolygon(3.25,.5)(2.75,0)(2.75,-.5)(3.25,0)
\rput[tl](3,-.2){$\frac{\sin \beta}{\sin \alpha}|d\vy|$}
\rput[l](3.35,.25){$\frac{\sin \beta-\alpha}{\sin \alpha}|d\vx|$}
\qdisk(.5,0){.04}
\qdisk(3,0){.04}

\endpspicture
\end{fullfigure}

The expectation \eqref{e:expect} can be computed explicitly.
Define angles $\alpha$ and $\beta$ by
$$\cos \alpha = \vx \cdot \vy \qquad \cos \beta = \vx \cdot \vp.$$ The cone
$C(d\vx,\vp)$ is an infinitesimally thin triangle which intersects $d\vy$ with
some probability.  If there is an intersection, its sign does not
depend $\vp$, but only on the sign of $\vx \wedge \vy \wedge d\vx \wedge d\vy.$
There is no intersection when $\beta < \alpha$.   For fixed $\beta > \alpha$,
the probability of an intersection, times its sign, is given by the
formula
$$\frac{(\sin \beta-\alpha)(\sin \beta)}{2\pi^2(\sin \alpha)^3}
    \vx \wedge \vy \wedge d\vx \wedge d\vy \wedge d\beta.$$
This can be seen by supposing that $d\vx$ and $d\vy$ are orthogonal to each
other and to both $\vx$ and $\vy$. Then the intersection window for $\vp$ is a
rectangle which is
$$\frac{\sin \beta-\alpha}{\sin \alpha}|d\vx| \qquad \mathrm{by} \qquad
\frac{\sin \beta}{\sin \alpha}|d\vy|.$$
(See \fig{f:triangle}.)  The factor $\vx \wedge \vy$ also produces a
numerical factor of $\sin \alpha$ which must be cancelled.  If we include
variation of $\beta$, the region for $\vp$ is a brick of thickness $d\beta$.
To get an overall expectation, we must integrate the volume of this
brick with respect to $\beta$ and divide by $\Vol S^3 = 2\pi^2$.  Using
$$\int_\alpha^\pi (\sin \beta-\alpha)(\sin \beta) d\beta
    = \frac{(\pi-\alpha)(\cos \alpha) + (\sin \alpha)}{2},$$
the final answer is
$$\omega = \frac{(\pi-\alpha)(\cos \alpha) + (\sin \alpha)}
    {4\pi^2(\sin \alpha)^3}\vx \wedge \vy \wedge d\vx \wedge d\vy.$$
This formula agrees with equation~\eqref{e:s3}.

\section{Open problems}

\subsection{Odds and ends}

If the convex body $K$ or its dual $K^\circ$ is not positively curved, then
$K^+$ is a weakly spacelike section of $H^+$, and $K^-$ is a weakly timelike
section of $H^-$. In other words, both will have some null tangents.  In general
if $N^\pm$ are weak necks of $H^\pm$ in $\R^{(a,b)}$, their join $N^+ * N^-$
might no longer be embedded except on a measure $0$.  For example, $K^+ * K^-$
is not embedded if $K$ is a polytope.

Nonetheless, formula~\eqref{e:int} is still the volume of a region
$N^\diamond$.  The necks are Lipschitz by virtue of being weakly spacelike and
timelike, so the integrand is $L^\infty$.  Moreover, the weak necks can be
approximated by strictly spacelike and timelike necks, so the integrand is never
negative; no volume is ever subtracted.  We conjecture that in this generality,
flat necks uniquely minimize the volume of $N^\diamond$.  This is immediate
from the  proof of \thm{th:abnecks} for strict necks, but weak necks introduce
new complications.

We do not know whether $K^\diamond$ is always convex when $K$ is a centrally
symmetric convex body.   In light of the Mahler conjecture, we conjecture that
$K^\diamond$ has the most volume when $K = C_n$.  (We do not even know this when
$n=2$.) A related phenomenon is the fact that $K^\diamond = K \times K^\circ$
when $K$ is a Hanner polytope \cite{Hanner:translates}.  We conjecture that they
are the only symmetric convex bodies with this property, in light of the
conjecture that they are the only convex bodies with minimum Mahler volume.

For two general necks $N^\pm$, $N^\diamond$ need not be convex and its volume is
unbounded even in fixed dimensions.  An upper bound may be possible if we
suppose that $N^\diamond$ lies between $H^+$ and $H^-$.  When $a = b = n$, then
$H^+$ and $H^-$ are contact manifolds, and the necks that come from convex
bodies are Legendrian.  We conjecture the converse (and it may not be hard to
prove), that every Legendrian neck of $H^+$ is $K^+$ for a convex body $K$. 
This extra geometric property of $K^+$ could be important for other problems in
convex geometry.

Reference~\citealp{Kuperberg:bottleneck} discusses a different generalized
bottleneck conjecture which is still open when $a,b > 2$.  Instead of minimizing
$$\int_{N^+ \times N^-} \vx \wedge \vy \wedge
    d\vx^{\wedge a-1} \wedge d\vy^{\wedge b-1},$$
we can minimize
$$\int_{N^+_1 \times N^+_2} (\vx \wedge d\vx^{\wedge a-1})
    \cdot (\vy \wedge d\vy^{\wedge a-1})$$
for two necks of $H^+$.  The dot product on $\Lambda^a(\R^{(a,b)})$ in this
integral is induced from the one on $\R^{(a,b)}$.  As before, we conjecture that
the minimum occurs when the necks are flat and centered, and in this case
coincident rather than orthogonal.  When $N^+ = K^+$ for a
symmetric convex body $K$, the two integrals are equal.

We can consider variations of the Mahler conjecture for complex and quaternionic
convex bodies, \ie, the unit balls of  complex and quaternionic norms.  The
natural conjecture is that the unit balls of the complex and quaternionic
$\ell^1$ and $\ell^\infty$ norms have the least Mahler volume.  (The 
$\ell^\infty$ ball is a polydisk, a Cartesian power of the Euclidean unit ball
or disk in the complex numbers or quaternions.)  Hanner polytopes also
generalize analogously to Hanner bodies.  If $K$ is a complex or quaternionic
convex body, then we can define not  only the necks $K^\pm$, but more generally
the neck $K^z$ for any scalar $|z| = 1$.  This extra geometry could lead to a
variant of \thm{th:abnecks} and an improvement to \cor{c:bound}.

\subsection{Philosophy}

We can attempt to apply the philosophy behind \thm{th:abnecks} and \cor{c:bound}
to other problems in asymptotic convex geometry.  Suppose that we want to bound
some affinely invariant function $f(K)$ on convex bodies.  Typically $f(K)$ is
maximized (respectively, minimized) by ellipsoids; the question is to find a
lower bound (respectively, an upper bound) to say that all convex bodies $K$ are
``not too far from round''.  Our philosophy is to find another roundness
statistic $g(K) \le f(K)$ (respectively $g(K) \ge f(K)$) which is
\emph{minimized} (respectively \emph{maximized}) when $K$ is an ellipsoid.  This
yields a bound for the original statistic $f(K)$.

For example, the isotropic constant conjecture asserts that if $K \subset \R^n$
is a convex body, its isotropic constant $L(K)$ is bounded above by some
constant $C$, independent of dimension \cite{Bourgain:maximal}. (See also
Reference \citealp{Ball:thesis}).  The centrally symmetric case of the
conjecture implies the general case \cite{Klartag:isom}. The constant can be
defined by minimizing an expectation:
$$L(K)^2 \stackrel{\mathrm{def}}{=}
    \min_{K' = T(K)} \frac{E_{K'}[\vx \cdot \vx]}{n(\Vol K')^{2/n}},$$
where the minimum is taken over all affine positions $K' = T(K)$, and the
expectation is taken with respect to a random point $\vx \in K'$. It is also not
hard to show that $L(K)$ is minimized when $K$ is an ellipsoid.  (This result is
due to Blaschke \cite{Blaschke:affine}, but the argument is very simple.  We can
minimize $E_K[\vx \cdot \vx]$ over all $K$ with fixed volume just by moving the
measure of $K$ as close to the origin as possible.)

\begin{conjecture} If $K \subset \R^n$ is centrally symmetric, then
$$E_{K \times K^\circ}[(\vx \cdot \vy)^2]$$
is maximized when $K$ is an ellipsoid.
\label{c:xy} \end{conjecture}

A number of convex geometers, including Ball and Giannopolous, have long known
that \conj{c:xy} implies the isotropic constant conjecture.  It is not hard to
show \cite[Prop.~1.2.3]{Giannopoulos:notes} that
$$E_{K \times K^\circ}[(\vx \cdot \vy)^2] \ge n^2 L(K)^2L(K^\circ)^2
(\Vol K)^{2/n}(\Vol K^\circ)^{2/n}.$$
The centrally symmetric isotropic constant conjecture then follows by combining
\conj{c:xy}, \thm{th:bm}, and Blaschke's inequality applied to $L(K^\circ)$.  On
the negative side, \conj{c:xy} also sharpens the Blaschke-Santal\'o theorem, so
it is more difficult and perhaps less likely.

We can propose even stronger and less likely conjectures than \conj{c:xy}.  We
conjecture that for every $0 < c < 1$, the probability
$$p_K(c) = P_{K \times K^\circ}[\vx \cdot \vy \ge c]$$
is maximized when $K$ is an ellipsoid.  We can even conjecture that
$$q_K(c) = P_{K \times K^\circ}[\vx \cdot \vy \ge
    c||\vx||_K \;||\vy||_{K^\circ}]$$
is maximized when $K$ is an ellipsoid; this would imply the conjecture for
$p_K(c)$.  It can be shown, using the technique of osculating ellipsoids, that
$$q_K(1-\eps) = \frac{(\Vol K^+)(\Vol B_{n-1})\eps^{\frac{n-1}2}}{v(K)}
    + o(\eps^{\frac{n-1}2}),$$
and that $p_K(1-\eps)$ has similar asymptotics.  Here $\Vol K^+$ is the volume
of $K^+$ as a spacelike, and therefore Riemannian, submanifold of $\R^n \times
\R^n$ with the inner product~\ref{e:ip}.  So we conjecture that
ellipsoids maximize the ratio $(\Vol K^+)/v(K)$, even though they maximize the
denominator. A related conjecture is the following:

\begin{conjecture} If $N^+ \subset H^+ \subset \R^{(a,b)}$ is a neck, then its
Riemannian volume $\Vol N^+$ is maximized when $N^+$ is flat.
\label{c:vol}
\end{conjecture}

We do not know if \conj{c:vol} is hard to prove, or even genuinely an
open problem, in the context of current methods in differential geometry
\cite{LS:slice}.

Note that if $N^+$ is only a weak neck, then $\Vol N^+$ can vanish.
For example, if $K$ is a polytope, then $\Vol K^+ = 0$.


\providecommand{\bysame}{\leavevmode\hbox to3em{\hrulefill}\thinspace}
\providecommand{\MR}{\relax\ifhmode\unskip\space\fi MR }
\providecommand{\MRhref}[2]{%
  \href{http://www.ams.org/mathscinet-getitem?mr=#1}{#2}
}
\providecommand{\href}[2]{#2}

\end{document}